\documentclass[twoside,11pt]{article}
\usepackage{pgm}

\usepackage{amsmath, amssymb}

\usepackage[utf8]{inputenc}
\inputencoding{utf8}

\usepackage{hyperref,color,soul,booktabs}
\setulcolor{blue}
\newcommand{\BibTeX}{\textsc{B\kern-0.1emi\kern-0.017emb}\kern-0.15em\TeX}

\newcommand{\sbsbsection}[1]{\noindent \textbf{#1.}  }

\ShortHeadings{The Variance of Causal Effect Estimators for Binary V-structures}{Kuipers and Moffa}

\begin{document}

\title{The Variance of Causal Effect Estimators for Binary V-structures}

\author{\Name{Jack Kuipers} \Email{jack.kuipers@bsse.ethz.ch} \\
\addr D-BSSE, ETH Zurich, Mattenstrasse 26, 4058 Basel, Switzerland
\and
\Name{Giusi Moffa} \Email{giusi.moffa@unibas.ch} \\
\addr Department of Mathematics and Computer Science, University of Basel, Basel, Switzerland \\
Division of Psychiatry, University College London, London, UK}

\maketitle

\begin{abstract}
Adjusting for covariates is a well established method to estimate the total causal effect of an exposure variable on an outcome of interest. Depending on the causal structure of the mechanism under study there may be different adjustment sets, equally valid from a theoretical perspective, leading to identical causal effects. However, in practice, with finite data, estimators built on different sets may display different precision. To investigate the extent of this variability we consider the simplest non-trivial non-linear model of a v-structure on three nodes for binary data.  We explicitly compute and compare the variance of the two possible different causal estimators.  Further, by going beyond leading order asymptotics we show that there are parameter regimes where the set with the asymptotically optimal variance does depend on the edge coefficients, a result which is not captured by the recent leading order developments for general causal models. As a practical consequence, the adjustment set selection needs to account for the relative magnitude of the relationships between variables with respect to the sample size, and cannot rely on purely graphical criteria.
\end{abstract}

\section{Introduction}

As graphical representations of multivariate probability distributions, Bayesian networks are widely used statistical models with an underlying directed acyclic graph (DAG) structure. When taking DAGs to represent causal diagrams \citep{greenland99causal, pearl00, hr06instruments, wr07four} we may use a machinery based on the `do' calculus of \cite{pearl95causal} to estimate potential intervention effects of any variable on any other.  Different graphical criteria exist to identify valid adjustment sets, among which the back-door criterion \citep{pearl93bayesian} is probably the most well known, and with more generalised strategies developed more recently \citep{svr10, perkovic17complete}. 

A valid adjustment set $\boldsymbol{Z}$ for the effect of $X$ on $Y$ is such that for any probability distribution $p$ compatible with the underlying graphical structure the probability distribution of $Y$ after intervening on $X$ (setting it to some value) satisfies \citep{mc15}
\begin{align}\label{eq:adjustment}
p(Y \mid \text{do } X) = \left\{ \begin{array}{lcr}
p(Y \mid X) && \text{if } \boldsymbol{Z} = \emptyset \\ 
\int_{\boldsymbol{z}}p(Y \mid X, \boldsymbol{z})p(\boldsymbol{z}) \mathrm{d}\boldsymbol{z}  & & \text{otherwise} 
\end{array}\right.
\end{align}
For linear Gaussian models, the marginalisation can be simply estimated by regressing $Y$ on $X$ and $\boldsymbol{Z}$ and extracting the coefficient of $X$, hence the naming of `adjustment' sets. This also holds for linear non-Gaussian causal models \citep{henckel19graphical}.

The set of parents of $X$ always satisfies the back-door criterion and is therefore a valid adjustment set, but there may be many more depending on the graphical structure of the DAG \citep{perkovic17complete}.  Although all valid adjustment sets provide consistent estimators of the causal effects, for finite-sized data different adjustment sets can lead to different numerical estimates, and with different precisions.

In evaluating the variance of different estimators, \cite{henckel19graphical} recently obtained the remarkable result that the asymptotically optimal adjustment set can be determined solely based on graphical criteria regardless of the edge coefficients.  Even more recently, this has been extended to non-parametric estimators \citep{rotnitzky2019efficient} and the asymptotically optimal set has been further characterised \citep{witte2020efficient}.

\begin{figure}[t]
 \centering
 \includegraphics[width=0.4\textwidth]{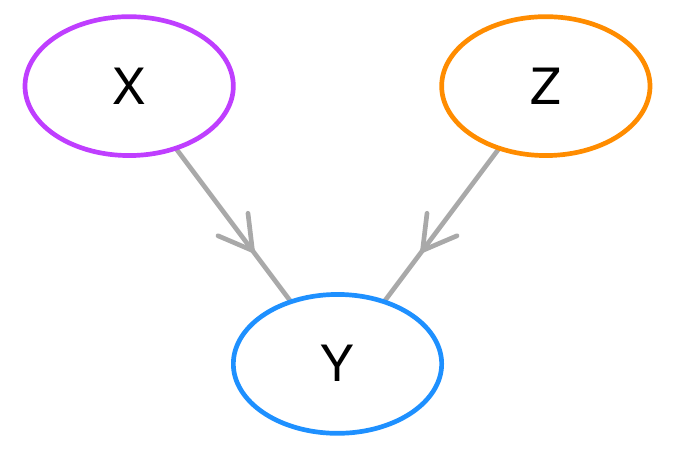}
\caption{A V-structure on 3 Nodes.}
\label{fig:vstruct}
\end{figure}

To explore the precision of causal estimators for non-linear models we consider the simplest such case: a DAG with 3 nodes of binary variables organised in a v-structure with the outcome $Y$ of interest as a collider with parents $Z$ and $X$ (Figure \ref{fig:vstruct}), and with the latter being the exposure whose effect we wish to estimate. For binary data and relatively small networks, one can explicitly marginalise over the remaining nodes in the DAG and its parameters \citep{moffaetal17} to derive interventional distributions as
\begin{align}
p(Y \mid \mathrm{do} (X) ) = \sum_{Z}p(Y, Z \mid \mathrm{do} (X))
\end{align}
and estimate causal effects from them.  

In the simple case of a v-structure (as in Figure \ref{fig:vstruct}) there is no confounding of the effect of $X$ on $Y$ (there are no common parents) so that the empty set constitutes a valid (and minimal) adjustment set, and the interventional distribution is simply
\begin{align}
p(Y \mid \mathrm{do} (X) ) = p(Y \mid X )
\end{align}
A valid expression for computing the total causal effect of $X$ on $Y$, in accordance with Equation (\ref{eq:adjustment}), is then
\begin{align}
F_{R} = p(Y \mid \mathrm{do}(X = 1)) - p(Y \mid \mathrm{do}(X = 0)) = p(Y \mid X = 1) - p(Y \mid X = 0)
\end{align}
where we used the subscript $R$ for raw, to highlight the fact that the formula only involves raw (or observed) conditional probabilities of $Y$ on $X$, which in this simple scenario are sufficient to identify the desired causal effects.

However, by definition, the conditional distribution of $Y$ on $X$ is also
\begin{align}
p(Y \mid X) = \sum_{Z}p(Y, Z \mid X) = \sum_{Z}p(Y \mid X, Z)p(Z)
\end{align}
Therefore another valid expression for the total causal effect of $X$ on $Y$ is
\begin{align}\label{eq:marginalised}
F_{M} = & P(Y \mid X = 1, Z = 1)P(Z) + P(Y \mid X = 1, Z = 0)(1-P(Z)) \nonumber \\
& - P(Y \mid X = 0, Z = 1)P(Z) - P(Y \mid X = 0, Z = 0)(1 - P(Z))
\end{align}
where we used the subscript $M$ to highlight the fact that the formula derives from explicitly marginalising $Z$ out from the joint distribution $p(Y, Z \mid X)$. In contrast, one could interpret the formula based on raw conditionals as performing the marginalisation implicitly (with the observations already providing a marginalised sample).

A more general way of understanding Equation (\ref{eq:marginalised}) is by observing that in the case of the v-structure $Z$ also constitutes a valid adjustment set (albeit not a minimal one). Then starting from the joint interventional distribution $p(Y, Z \mid \mathrm{do} (X) )$ the interventional distribution of $Y$ when intervening on $X$ is also
\begin{align}
p(Y \mid \mathrm{do} (X) ) = \sum_{Z}p(Y \mid \mathrm{do} (X), Z)p(Z \mid \mathrm{do} (X)) = \sum_{Z}p(Y \mid X, Z)p(Z)
\end{align}
with the latter equality justified by structural and invariance properties, and also in agreement with the standardisation formula in Equation (\ref{eq:adjustment}).

Since we see that adjustment by $Z$ is valid, but not necessary, it is natural to ask whether the two estimators differ in terms of precision. To answer the question we compute the variance, for finite sample sizes, of the two different estimators corresponding to the implicit and explicit marginalisation as outlined before.

It is instructive to also consider the DAG with the edge from $Z\to Y$ deleted.  Since $Y$ would then be independent of $Z$, the marginalisation would reduce to the raw conditionals. The estimator using raw conditionals is therefore the same whether the edge from $Z\to Y$ is present or not, while the approach using marginalisation would give different estimates for the two cases. Intuitively we would expect that the extent by which estimates differ will depend on the strength of the relationship between $Y$ and $Z$. The underlying rationale is similar to that for the standard practice of including baseline covariates in linear models of the outcome in randomised controlled trials \citep{senn2011modelling}, where prognostic factors and the (randomised) treatment can be seen as forming a v-structure with the outcome as the collider.

The v-structure actually provides the simplest example where there is a choice between different adjustment sets. If we add an edge in the graph of Figure \ref{fig:vstruct} connecting $X$ and $Z$ we end up with no choice about adjustment sets: in particular if we add an edge from $X \to Z$, then $Z$ is not a valid adjustment set and the empty set is the only choice; conversely if we add an edge from $Z \to X$ then $Z$ is a confounder (a common parent) and it must be adjusted for, making it the only valid adjustment set with the empty set no longer valid.

\section{Causal Estimates for a Binary V-structure}

For both causal estimators we will use the maximum likelihood estimates of probabilities from the observed data. We consider the DAG in Figure \ref{fig:vstruct} with the following probability tables:
\begin{align}
p(X = 1) = p_X \, , & & p(Y = 1 \mid X = 0, Z = 0) &= p_{Y, 0}  \, , & & p(Y = 1 \mid X = 1, Z = 0) = p_{Y, 2} \nonumber \\
p(Z = 1) = p_Z \, , & & p(Y = 1 \mid X = 0, Z = 1) &= p_{Y, 1}  \, , & & p(Y = 1 \mid X = 1, Z = 1) = p_{Y, 3}
\end{align}

When we generate data, as a collection of $N$ binary vectors, from the DAG in Figure \ref{fig:vstruct}, instead of forward sampling along the topological order for this small example we can sample directly from a multinomial with probabilities
\begin{equation}
\begin{array}{c|c|c|lcc|c|c|l}
X & Z & Y & p & & X & Z & Y & p \\
\cline{1-4}
\cline{6-9}
0 & 0 & 0 & p_0 = (1-p_X)(1-p_Z)(1-p_{Y,0}) & & 1 & 0 & 0 & p_4 = p_X(1-p_Z)(1-p_{Y,2}) \\
0 & 0 & 1 & p_1 = (1-p_X)(1-p_Z)p_{Y,0} & & 1 & 0 & 1 & p_5 = p_X(1-p_Z)p_{Y,2} \\
0 & 1 & 0 & p_2 = (1-p_X)p_Z(1-p_{Y,1}) & & 1 & 1 & 0 & p_6 = p_Xp_Z(1-p_{Y,3}) \\
0 & 1 & 1 & p_3 = (1-p_X)p_Zp_{Y,1} & & 1 & 1 & 1 & p_7 = p_Xp_Zp_{Y,3}
\end{array}
\end{equation}
If we represent with $N_i$ the number of sampled binary vectors indexed by $i = 4X+2Z+Y$, then the estimator of $F$ from the raw conditionals is simply
\begin{align}
R = R_1 - R_0\, , & & R_1 = \frac{N_5 + N_7}{N_4 + N_5 + N_6 + N_7}\, , & & R_0 = \frac{N_1 + N_3}{N_0 + N_1 + N_2 + N_3}
\end{align}

Using the marginalisation we would have the following estimator
\begin{align}
M = M_1 - M_0 \, , & &
M_1 = M_{11} + M_{10} \, , & &  M_0 = M_{01} + M_{00}
\end{align}
with the terms separated for later ease
\begin{align}
M_{11} = \frac{N_7}{(N_6+N_7)}\frac{(N_2 + N_3 + N_6 + N_7)}{N} \, , & & M_{01} = \frac{N_3}{(N_2+N_3)}\frac{(N_2 + N_3 + N_6 + N_7)}{N} \nonumber \\
M_{10} = \frac{N_5}{(N_4+N_5)}\frac{(N_0 + N_1 + N_4 + N_5)}{N} \, , & & M_{00} = \frac{N_1}{(N_0+N_1)}\frac{(N_0 + N_1 + N_4 + N_5)}{N}
\end{align}

These estimators, as they rely on observed data frequencies, are non-parametric and fit in the general framework of \cite{rotnitzky2019efficient} for arbitrary graphs.
The key advance of our derivation with respect to their result is that we consider terms beyond the leading order asymptotics and compute
the variance of the estimators for arbitrary sample sizes, which further enables us to perform more detailed asymptotic analyses.

\subsection{Raw Conditionals}

To compute $E[R]$ we need to average over a multinomial sample
\begin{equation}
E[R] = \sum \frac{N!}{N_0! \cdots N_7!}p_0^{N_0}\cdots p_7^{N_7} R
\end{equation}
for which we use that fact that $(p_0 + \ldots + p_7)^N$ generates the probability distribution when we perform a multinomial expansion. To obtain the terms needed for the expectation we define
\begin{align} \label{SNcomp}
S_N = & \left\{[p_0+p_2 + (p_1+p_3)w]x + [p_4+p_6 + (p_5+p_7)v]z\right\}^{N}
\end{align}
whose expansion is
\begin{align} \label{SNmulti}
S_N = & \sum \frac{N!}{N_0! \cdots N_7!}p_0^{N_0}\cdots p_7^{N_7} w^{N_1 + N_3}x^{N_0 + N_1 + N_2 + N_3} v^{N_5 + N_7}z^{N_4 + N_5 + N_6 + N_7}
\end{align}
Setting all the generating variables to 1 removes them from consideration and the generating function simplifies to the value 1:
\begin{align}
S_N\Big|_{\substack{w=x=1 \cr v=z=1}} = 1
\end{align}

The advantage of using generating functions \citep{wilf2005} is that we can express expectations in terms of differential and integral operators. For example the operator $v \frac{\partial}{\partial v}$ acting on $S_N$ will bring down a factor of $(N_5 + N_7)$ from the power of $v$ through differentiation and we then multiply by $v$ to leave the power unchanged (a useful feature to apply multiple operators later).  The effect of the operator is easiest to see when we apply it to the expanded form of $S_N$ from Equation (\ref{SNmulti}):
\begin{align}
v \frac{\partial}{\partial v} S_N = & \sum \frac{N!}{N_0! \cdots N_7!}p_0^{N_0}\cdots p_7^{N_7} (N_5 + N_7) w^{N_1 + N_3}x^{N_0 + N_1 + N_2 + N_3} v^{N_5 + N_7}z^{N_4 + N_5 + N_6 + N_7}
\end{align}
Removing the generating variables after applying the operator leads to
\begin{align}
v \frac{\partial}{\partial v} S_N \Big|_{\substack{w=x=1 \cr v=z=1}} = & \sum \frac{N!}{N_0! \cdots N_7!}p_0^{N_0}\cdots p_7^{N_7} (N_5 + N_7) = E[N_5 + N_7]
\end{align}
which is an expectation over the multinomial probability distribution of a binary v-structure. To actually perform the differentiation we employ the compact form of $S_N$ from Equation (\ref{SNcomp}) to easily obtain the result of $N(p_5 + p_7)$.

The integral operator $\int \mathrm{d} z \frac{1}{z}$ will introduce a denominator of $(N_4 + N_5 + N_6 + N_7)$, so by combining operators we derive our first expectation of interest
\begin{align}
E[R_1] = \int \frac{v}{z}\frac{\partial}{\partial v} S_N \, \mathrm{d}z \, \Bigg|_{\substack{w=x=1 \cr v=z=1}} = \frac{v(p_5+p_7)}{p_4+p_6+(p_5+p_7)v}S_N \, \Bigg|_{\substack{w=x=1 \cr v=z=1}}
\end{align}
When we substitute for the generating variables (which sets $S_N=1$), and perform the same steps for $R_0$ we obtain
\begin{align}
E[R] = \frac{p_5+p_7}{p_4+p_5+p_6+p_7} - \frac{p_1+p_3}{p_0+p_1+p_2+p_3} = \frac{p_5+p_7}{p_X} - \frac{p_1+p_3}{1-p_X}
\end{align}

\sbsbsection{The Variance}
To compute the variance 
\begin{equation}
V[R] = V[R_1] - 2C[R_1, R_0] + V[R_0]
\end{equation}
we first show that the covariance is 0
\begin{align}
E[R_1R_0] &= \int \frac{w}{x} \frac{\partial}{\partial w} \int \frac{v}{z}\frac{\partial}{\partial v} S_N \, \mathrm{d}z \, \mathrm{d}x \, \Bigg|_{\substack{w=x=1 \cr v=z=1}} = \int \frac{w}{x} \frac{\partial}{\partial w} \frac{v(p_5+p_7)}{p_4+p_6+(p_5+p_7)v}S_N \, \mathrm{d}x \, \Bigg|_{\substack{w=x=1 \cr v=z=1}} \nonumber \\
&= \frac{v(p_5+p_7)}{p_4+p_6+(p_5+p_7)v} \cdot \frac{w(p_1+p_3)}{p_0+p_2+(p_1+p_3)w}S_N \, \Bigg|_{\substack{w=x=1 \cr v=z=1}} = E[R_1]E[R_0] 
\end{align}
The last equality follows by comparing to the values of $E[R_1]$ and $E[R_0]$ computed above.

The more tricky terms are
\begin{align}
E[R_1^2] &= \int \frac{v}{z} \frac{\partial}{\partial v} \int \frac{v}{z}\frac{\partial}{\partial v} S_N \, \mathrm{d}z \, \mathrm{d}z \, \Bigg|_{\substack{w=x=1 \cr v=z=1}} = \int \frac{v}{z} \frac{\partial}{\partial v} \frac{v(p_5+p_7)}{p_4+p_6+(p_5+p_7)v}S_N \, \mathrm{d}z \, \Bigg|_{\substack{w=x=1 \cr v=z=1}} \nonumber \\
&= \frac{v(p_5+p_7)}{p_4+p_6+(p_5+p_7)v}\int \left[ \frac{(p_4+p_6)}{p_4+p_6+(p_5+p_7)v}\frac{S_N}{z} + \frac{v}{z} \frac{\partial}{\partial v} S_N \right] \, \mathrm{d}z \, \Bigg|_{\substack{w=x=1 \cr v=z=1}} \nonumber \\
&= \frac{(p_5+p_7)(p_4+p_6)}{p_{X}^2}\int \frac{\left(1-p_{X} + p_{X}z\right)^{N}}{z} \, \mathrm{d}z \, \Bigg|_{{z=1}} + E[R_1]^2
\end{align}
The remaining integral can be expressed in terms of hypergeometric functions:
\begin{align} \label{eq:hypergeom}
\int \frac{\left(1-p_{X} + p_{X}z\right)^{N}}{z} \, \mathrm{d}z &= \sum_{k=1}^{N} \binom{N}{k}\frac{1}{k}p_X^k(1-p_X)^{N-k} \nonumber \\
&= Np_X(1-p_X)^{N-1}F\left([1, 1, 1-N], [2, 2], -\frac{p_X}{1-p_X}\right)
\end{align}
where we use the notation $F\left([a_1,\ldots, a_p], [b_1, \ldots, b_q], t\right)$ for the generalised hypergeometric function ${}_{p}F_q\left(a_1,\ldots, a_p; b_1, \ldots, b_q; t\right)$ with square brackets to help delineate the different arguments.
Repeating the calculations for $V[R_0]$ we obtain
\begin{align}
V[R] = & \frac{(p_5+p_7)(p_4+p_6)}{p_{X}}N(1-p_X)^{N-1}F\left([1, 1, 1-N], [2, 2], -\frac{p_X}{1-p_X}\right) \nonumber \\
& + \frac{(p_1+p_3)(p_0+p_2)}{(1-p_{X})}Np_X^{N-1}F\left([1, 1, 1-N], [2, 2], -\frac{1-p_X}{p_X}\right)
\end{align}
We discuss bounds on this variance in Appendix \ref{Rbounds}.

\subsection{Marginalisation}

To compute the expected value $E[M]$  we define
\begin{align}
T_N = & \left\{ap_0s + ap_1st + bp_2u + bp_3uv + ap_4w + ap_5wx + bp_6y + bp_7yz\right\}^{N}
\end{align}
where we include extra generating variables for all terms in our estimators. Then
\begin{align} \label{Mexpected}
E[M_{11}] &= \frac{1}{N}\left[\int\mathrm{d}y \frac{bz}{y}\frac{\partial^2}{\partial b\partial z} \right] T_N \, \Bigg|_{\substack{a=b=1 \cr s=t=1 \cr u=v=1 \cr w=x=1 \cr y=z=1}} = \left[\frac{bzp_7(p_2u+p_3uv + p_6y + p_7yz)}{(p_6 + p_7z)}\right] T_{N-1} \, \Bigg|_{\substack{a=b=1 \cr s=t=1 \cr u=v=1 \cr w=x=1 \cr y=z=1}} \nonumber \\
&= \frac{p_7}{(p_6 + p_7)}p_Z
\end{align}
and similarly for the other terms, leading to
\begin{align}
E[M] = & \frac{p_7}{(p_6 + p_7)}p_Z + \frac{p_5}{(p_4 + p_5)}(1-p_Z) - \frac{p_3}{(p_2 + p_3)}p_Z - \frac{p_1}{(p_0 + p_1)}(1-p_Z)
\end{align}

To compute the variance, we reapply the operators of Equation (\ref{Mexpected}), as detailed in Appendix \ref{Mvariance}.

\begin{figure}
 \centering
 \begin{tabular}{cc}
 \includegraphics[width=0.45\textwidth]{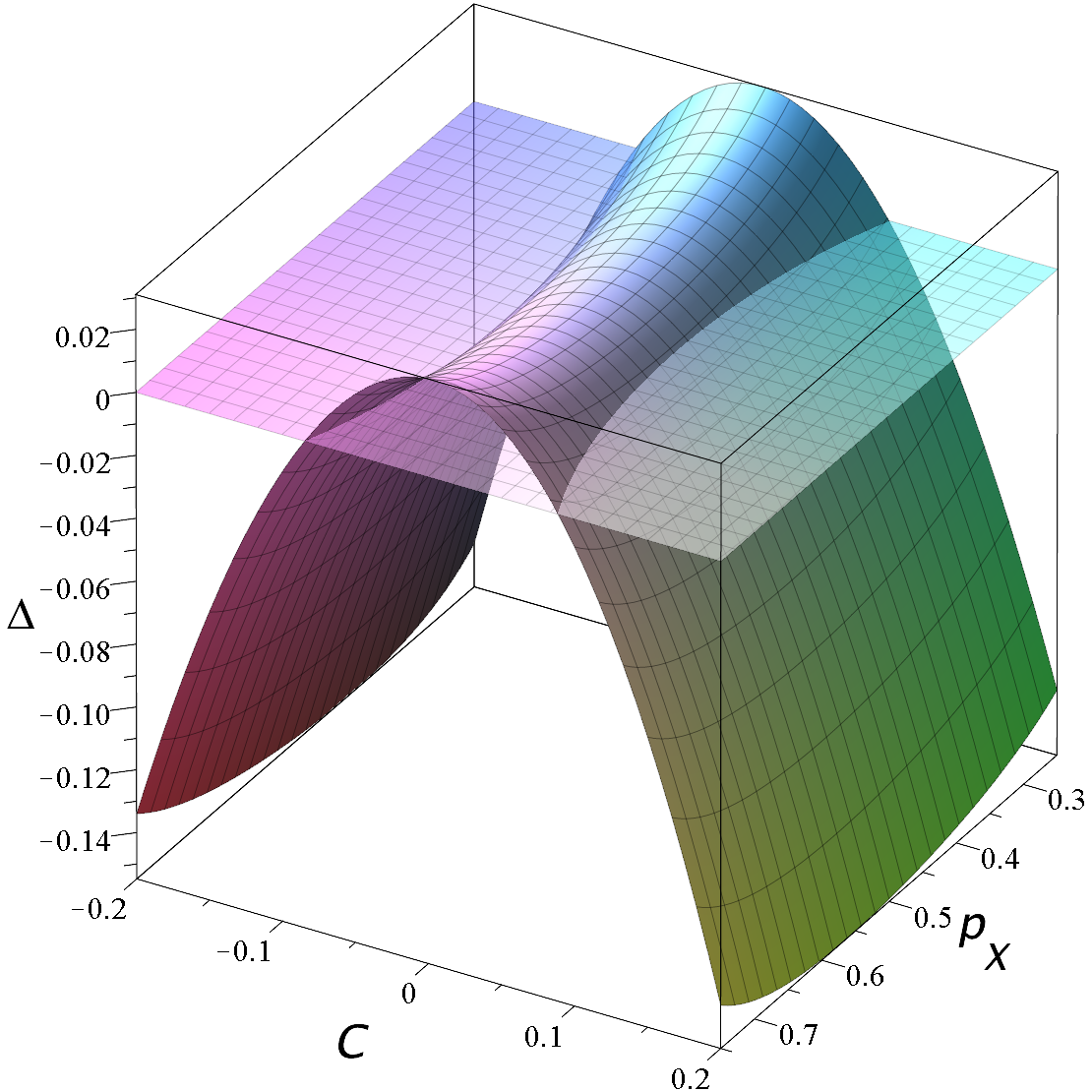} & \includegraphics[width=0.45\textwidth]{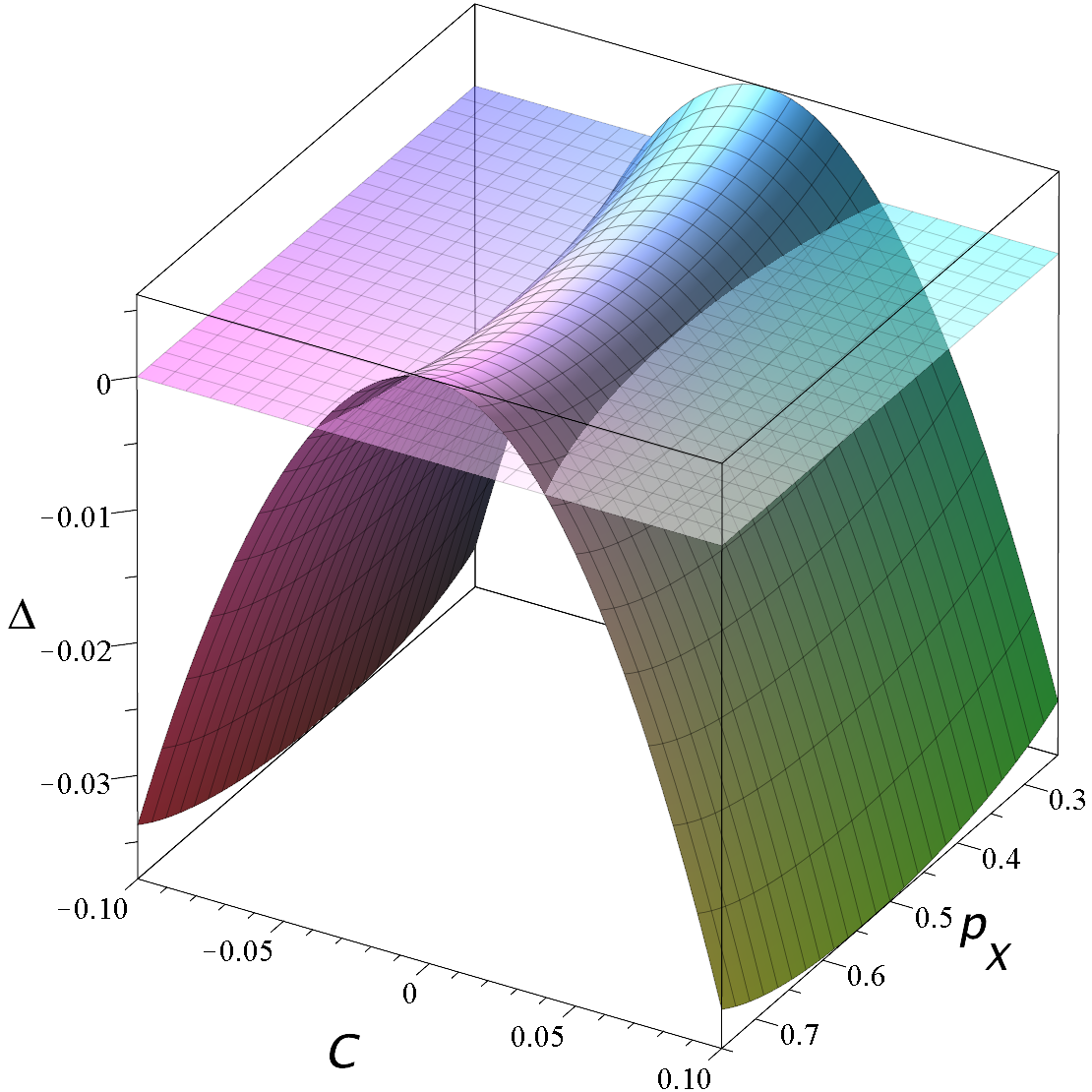} \\
 (a) $N=100$ & (b) $N=400$
 \end{tabular}
\caption{The Relative Difference in Variance of the Two Estimators.}
\label{fig:deltaplot}
\end{figure}

\subsection{Numerical checks}

Code to evaluate the variance of the two estimators through simulation, as well as to evaluate the analytical results above is hosted at \url{https://github.com/jackkuipers/Vcausal}. As an example, for $p_X = \frac{1}{3}, p_Z = \frac{2}{3}, p_{Y,0} = \frac{1}{6}, p_{Y,1} = \frac{1}{2}, p_{Y,2} = \frac{1}{3}, p_{Y,3} = \frac{5}{6}$ and $N=100$, we obtained Monte Carlo estimates of the standard deviation of $R$ and $M$ as $0.101929$ and $0.0924014$ respectively from 40 million repetitions. This agrees with the respective analytical results of $0.101932$ and $0.0924017$.

\subsection{Relative Difference in Variances}

To explore the difference in variances, we write the probabilities as
\begin{align}
p_{Y, 0} &= q_0 - C \, , & p_{Y, 1} &= q_0 + C \, , &
p_{Y, 2} &= q_1 - C \, , & p_{Y, 3} &= q_1 + C 
\end{align}
where $C$ is a measure of the effect of $Z$ on $Y$ (the same for each $X$) and the causal effect of $X$ on $Y$ is $q_1 - q_0$.

We plot the difference in variances of the two estimators, $\Delta = \frac{V[M]- V[R]}{V[R]}$.  In Figure \ref{fig:deltaplot} we leave $p_X$ free, set $p_Z=\frac{2}{3}$ and set $q_0=\frac{1}{3}$, $q_1=\frac{2}{3}$ and plot $\Delta$ for $N=100$ and $N=400$.  In the plot for $N=400$ we also scaled $C$ by dividing by 2. The behaviour and rescaled plots are very similar, suggesting a $N^{-\frac{1}{2}}$ scaling.

\section{Asymptotic Behaviour}

To examine the asymptotic behaviour of the causal effect estimators in more detail, we expand the hypergeometric function as in Appendix \ref{asymphyper}, to obtain the following for the variance of $R$
\begin{align}
V[R] \cdot N = & \frac{(q_1 + (2p_Z - 1)C)(1-q_1 - (2p_Z - 1)C)}{p_X}\left(1 + \frac{(1-p_X)}{Np_X}\right) \nonumber \\
& + \frac{(q_0 + (2p_Z - 1)C)(1-q_0 - (2p_Z - 1)C)}{(1-p_X)}\left(1 + \frac{p_X}{N(1-p_X)}\right) + O(N^{-2})
\end{align}
and for the variance of $M$
\begin{align}
V[M] \cdot N = &\frac{q_1(1-q_1) - C^2}{p_X}\left(1 + \frac{2(1-p_X)}{Np_X}\right) - \frac{(2q_1 - 1)(2p_Z - 1)C}{p_X} \nonumber \\
& + \frac{q_0(1-q_0) - C^2}{1-p_X}\left(1 + \frac{2p_X}{N(1-p_X)}\right) - \frac{(2q_0 - 1)(2p_Z - 1)C}{(1-p_X)} + O(N^{-2})
\end{align}
To extract the asymptotic behaviour of the difference in variances of the two estimators, we consider $C \sim N^{-\frac{1}{2}}$ to obtain
\begin{align}
\left(V[M] -V[R]\right)\cdot N = &\frac{q_1(1-q_1)(1-p_X)}{Np_X^2} + \frac{q_0(1-q_0)p_X}{N(1-p_X)^2} - \frac{4p_Z(1-p_Z)}{p_X(1-p_X)}C^2 + O(N^{-\frac{3}{2}})
\end{align}
with root
\begin{align}\label{cstar}
C^{*} = \sqrt{\frac{1}{4Np_Z(1-p_Z)}}\left[\frac{q_1(1-q_1)}{p_X} (1-p_X)^2+ \frac{q_0(1-q_0)}{(1-p_X)}p_X^2\right]^{\frac{1}{2}}
\end{align}
so that
\begin{align} \label{varcross}
\lim_{N \to \infty} V[M] - V[R] < 0 \, , & & C > C^{*} \nonumber \\
\lim_{N \to \infty} V[M] - V[R] > 0 \, , & & C < C^{*}
\end{align}

Note that although we used the scaling $C \sim N^{-\frac{1}{2}}$ to extract this result, it holds more generally.  For example for fixed $C\neq 0$, it is trivial to see that $C > C^{*}$ for some $N$ and so that $V[M]$ will become lower than $V[R]$ in the limit $N\to\infty$.  The asymptotically optimal adjustment set therefore uses marginalisation rather than raw conditioning, in line with previous results \citep{henckel19graphical, rotnitzky2019efficient} from the leading order asymptotics. For fixed $C=0$ however, raw conditioning would be better. It is exactly by treating subleading terms, as we do here, that we can examine where the transition occurs and how it depends on the coefficients.  For weaker effects of the edge from $Z \to Y$, with $C \lesssim N^{-\frac{1}{2}}$, the raw conditional can give a more precise estimate of the causal effect of $X$ on $Y$.

The result in Equation (\ref{varcross}) therefore shows that the optimal adjustment set, for the effect of $X$ on $Y$ does not depend solely on graphical criteria, but crucially on the relative effect size of the influence of the other node $Z$ on the outcome $Y$. The result details, for large but finite sample sizes, the crossover in efficiency from $R$ to $M$ as the relative effect of $C$ increases. Therefore, starting from plausible parameter values we can use Equations (\ref{cstar}) and (\ref{varcross}) to guide our choice of adjustment. 

Since for many practical purposes the sample size may be determined by other considerations and we can never take the limit $N\to\infty$ the asymptotic regime developed here accounting for the relative scale of effects compared to the sample size is most relevant. Though derived just for the binary v-structure, this is a counterexample showing that recent leading order aysmptotic results \citep{rotnitzky2019efficient} cannot directly extend outside their particular asymptotic limit.

\section{Implications for causal discovery}

With a larger sample size, we may be able to detect and quantify smaller causal effects.  Therefore we wish to get a feeling for the strength of the edge $Z \to Y$ we would detect from the data, or equivalently for which values of $C$ we would infer the presence of the edge.
To do so, we calculate the expected difference in maximised log-likelihoods when including the edge compared to a DAG with the edge deleted:
\begin{align}
E[\Delta l] = & \frac{1}{2} + N_7 \ln \left(q_1 + C\right) + N_6 \ln \left(1 - q_1 - C\right) - N_7 \ln \left(q_1\right) - N_6 \ln \left(1 - q_1\right) + \ldots \nonumber \\
= & \frac{1}{2} + Np_X p_Z\left(q_1 + C\right)\ln \left(1 + \frac{C}{q_1}\right) + Np_X p_Z\left(1 - q_1 - C\right)\ln \left(1 - \frac{C}{1 - q_1}\right) + \ldots \nonumber \\
= & \frac{1}{2} + \frac{N}{2}\left[\frac{p_X}{q_1(1-q_1)}+\frac{(1-p_X)}{q_0(1-q_0)}\right]C^2 + O(C^3)
\end{align}
where the $\frac{1}{2}$ comes from Wilk's theorem \citep{wilks38large} for the additional parameter when maximising all the probabilities relative to evaluating with the restriction $C=0$.

The change is AIC is then
\begin{align}
E[\Delta \mathrm{AIC}] = 2 - 2E[\Delta l] = 1 - N\left[\frac{p_X}{q_1(1-q_1)}+\frac{1-p_X}{q_0(1-q_0)}\right]C^2 + O(C^3)
\end{align}
There is therefore an asymptotic regime where the edge is strong enough to detect on average using the AIC but the estimator from raw conditionals that does not use the edge has lower variance
\begin{align}
N(C^{*})^{2} \geq NC^2 \geq \left[\frac{p_X}{q_1(1-q_1)}+\frac{(1-p_X)}{q_0(1-q_0)}\right]^{-1} 
\end{align}
which follows from the Cauchy-Schwarz inequality. The regime only vanishes when $p_Z = \frac{1}{2}$ and $q_1(1-q_1)(1-p_X)^2 = q_0(1-q_0)p_X^2$ and the two bounds become equal.  Utilising the BIC instead ($E[\Delta \mathrm{BIC}] = E[\Delta \mathrm{AIC}] + \log(N)$) leads to a large regime where we would not detect the edge on average, but where the estimator using marginalisation that does rely on the edge has lower variance.

\section{Discussion}

To evaluate the precision of different estimators targeting the same causal effect in causal diagrams, we considered the simple case of a v-structure for binary data and explicitly computed the variance of the two different estimators for the effect of $X$ on the collider $Y$, with $Z$ as the other parent.
 
The results involve combinations of hypergeometric functions, suggesting that exact results for larger DAGs may be rather complex.  Which estimator has the lower variance depends, among other parameters, on the relative strength of the edge from $Z$ to $Y$.  In general, estimating the causal effect through marginalisation offers better performance in the presence of a stronger direct effect of $Z$ on $Y$. When the direct effect is weaker instead, ignoring the edge and estimating the causal effect through the raw conditionals provides higher precision.

By examining the asymptotic regime of large sample sizes, we could confirm the intuition that for edge strengths statistically detectable by the AIC, accounting for the edge in the estimation should generally lead to lower variance.  Conversely, that the presence of statistically non-detectable edges should be ignored to achieve a lower variance.  

Most importantly, we could also discover an asymptotic regime where raw conditional estimates, ignoring the edge, were more precise in the presence of statistically detectable edges. One way to appreciate the practical relevance of these findings is by observing that we can expect ranges of causal strengths which become statistically detectable from data before we can gain precision by accounting for them in the estimation. Our detailed asymptotic analysis for the v-structure goes beyond the leading-order asymptotic result where the optimal estimator does not depend on the edge coefficients \citep{henckel19graphical, rotnitzky2019efficient}.

Outside the asymptotic regime, for finite sample sizes the gain in precision when using marginalisation and thus explicitly accounting for the edge presence, appears to be linked to its strength.  Although the example considered here is the simplest non-trivial DAG, this finding further supports the idea that learning the full structure of the graph, beyond simply identifying a valid adjustment set, may benefit the precision of causal inference.  The practical limitation with observational data is that we can only learn structures up to an equivalence class, so that we need to consider the possible range of causal effects across the whole class \citep{mkb09}, or implement Bayesian model averaging across DAGs \citep{moffaetal17}.

If we use a more stringent criterion to decide about the presence of edges, such as the BIC for example, which implements a stronger penalisation with respect to the AIC, we may end up missing edges too weak to detect on average, but whose presence would improve the precision of the causal estimation through marginalisation. In other words, for moderately weak direct effects, the selection of suitable adjustment sets may be relatively sensitive to the choice of the score.  Analogously, we may expect that optimal causal estimation may also be sensitive to the choice of learning algorithm, whether constraint-based \citep{bk:sgs00, art:KalischB2007}, score-based search \citep{chickering02} or Bayesian sampling \citep{fk03, km17, ksm18}.  Quantifying the extent by which the structure learning affects causal estimation constitutes an interesting line of further investigation.

\vskip 0.2in
\bibliography{Vcausal}

\appendix

\section{Bounds on the variance of the $R$ estimator} \label{Rbounds}

This hypergeometric function in Equation (\ref{eq:hypergeom}) has a maximum value at around $\frac{1.5}{N}$, and we note that if we divide by $(k+1)$ instead of $k$ in the sum we have the simple result
\begin{align}
\sum_{k=0}^{N} \binom{N}{k}\frac{1}{k+1}p_X^k(1-p_X)^{N-k} = \frac{1}{p_X(N+1)} - \frac{(1-p_X)^{N+1}}{p_X(N+1)} 
\end{align}
so that by considering the early terms in the sum we can bound
\begin{align}
\sum_{k=1}^{N} \binom{N}{k}\frac{1}{k}p_X^k(1-p_X)^{N-k} > \frac{1}{p_X(N+1)} \, , & & 
p_X > \frac{N-1 + \sqrt{3N^2 + 4N+1}}{N(N+3)}
\end{align}
which we can loosen to $p_X > \frac{1+\sqrt{3}}{N}$.  This provides the following lower bound for the variance
\begin{align}
&V[R] > \frac{(p_5+p_7)(p_4+p_6)}{p_{X}^3(N+1)}  + \frac{(p_1+p_3)(p_0+p_2)}{(1-p_{X})^{3}(N+1)} \, , & & \frac{1+\sqrt{3}}{N} < p_X <  \frac{N - 1 - \sqrt{3}}{N}
\end{align}

To obtain a simple upper bound we can compute
\begin{align}
\sum_{k=1}^{N} \binom{N}{k}\frac{1}{k}p_X^k(1-p_X)^{N-k} 
< 2\sum_{k=1}^{N} \binom{N}{k}\frac{1}{k+1}p_X^k(1-p_X)^{N-k} 
< \frac{2}{p_X(N+1)}
\end{align}
so that the variance vanishes in the large $N$ limit
\begin{align}
V[R] &< \frac{2(p_5+p_7)(p_4+p_6)}{p_{X}^3(N+1)}  + \frac{2(p_1+p_3)(p_0+p_2)}{(1-p_{X})^{3}(N+1)}
\end{align}

\section{The variance of the $M$ estimator} \label{Mvariance}

For computing the variance of $M$, we need to reapply the operators used to obtain the expected value as in Equation (\ref{Mexpected}). If they act on different generating variables, they will simply recreate terms like the mean, so we focus on terms where they repeat.

\sbsbsection{A Variance}
For example:
\begin{align}
E[M_{11}^2] \cdot N = \int\mathrm{d}y \frac{bz}{y}\frac{\partial^2}{\partial b\partial z} \left[\frac{bzp_7u(p_2+p_3v)}{(p_6 + p_7z)} + bzp_7y\right] T_{N-1} \, \Bigg|_{\substack{a=b=1 \cr s=t=1 \cr u=v=1 \cr w=x=1 \cr y=z=1}} 
\end{align}
For the linear term in $y$, it is easiest if we rearrange and integrate first
\begin{align}
bz\frac{\partial^2}{\partial z \partial b} \int \mathrm{d} y bzp_7 \, T_{N-1} \, \Bigg|_{\substack{a=b=1 \cr s=t=1 \cr u=v=1 \cr w=x=1 \cr y=z=1}} 
= \frac{p_6p_7}{(p_6 + p_7)^2}p_Z + \frac{p_7^2}{(p_6 + p_7)} + (N-1)\frac{p_7^2}{(p_6 + p_7)}p_{Z}
\end{align}
while for the rest of $E[M_{11}^2]$ we first differentiate wrt $b$
\begin{align}
b\frac{\partial}{\partial b}\left[\frac{bzp_7u(p_2+p_3v)}{(p_6 + p_7z)}\right] T_{N-1} \, \Bigg|_{\substack{a=b=1 \cr s=t=1 \cr u=v=1 \cr w=x=1}}
= & (N-1)\left[\frac{zp_7(p_2+p_3)^2}{(p_6 + p_7z)} + zp_7(p_2+p_3)y\right]T_{N-2} \nonumber \\
& {} + \left[\frac{zp_7(p_2+p_3)}{(p_6 + p_7z)}\right] T_{N-1} \nonumber \\
\end{align}
For the part with the factor of $y$, we again integrate first wrt $y$ and then differentiate to obtain
\begin{align}
z\frac{\partial}{\partial z} \int \mathrm{d} y (N-1)zp_7(p_2+p_3) \, T_{N-2} \, \Bigg|_{\substack{a=b=1 \cr s=t=1 \cr u=v=1 \cr w=x=1 \cr y=z=1}} 
=& (p_2+p_3)\left[\frac{p_6p_7}{(p_6 + p_7)^2} + (N-1)\frac{p_7^2}{(p_6 + p_7)}\right]
\end{align}
on the rest we apply the operator for $z$
\begin{align}
z\frac{\partial}{\partial z} \ldots \Bigg|_{z=1} = & \left[\frac{p_6p_7(p_2+p_3)}{(p_6 + p_7)^2}\right] T_{N-1} + (N-1)\left[\frac{p_7^2(p_2+p_3)}{(p_6 + p_7)}y + \frac{p_6p_7(p_2+p_3)^2}{(p_6 + p_7)^2}\right] T_{N-2}  \nonumber \\
& {} + (N-1)(N-2)y\frac{p_7^2(p_2+p_3)^2}{(p_6 + p_7)}T_{N-3}
\end{align}
The linear terms in $y$ give the following
\begin{align}
\frac{p_7^2}{(p_6+p_7)^2}(p_2+p_3) + (N-1)\frac{p_7^2}{(p_6+p_7)^2}(p_2+p_3)^2
\end{align}
while the integrals lead to
\begin{align}
& \frac{p_6p_7(p_2+p_3)}{(p_6+p_7)}(N-1)(1-p_6-p_7)^{N-2} F\left([1, 1, 2-N], [2, 2], -\frac{p_6+p_7}{1-p_6-p_7}\right) \nonumber \\ 
& + \frac{p_6p_7(p_2+p_3)^2}{(p_6+p_7)}(N-1)(N-2)(1-p_6-p_7)^{N-3} F\left([1, 1, 3-N], [2, 2], -\frac{p_6+p_7}{1-p_6-p_7}\right)
\end{align}
Combining all the terms, subtracting the mean part squared and simplifying slightly we obtain
\begin{align}
& V[M_{11}] \cdot N = \frac{p_6p_7(p_2+p_3)}{(p_6+p_7)}(N-1)(1-p_6-p_7)^{N-2} F\left([1, 1, 2-N], [2, 2], -\frac{p_6+p_7}{1-p_6-p_7}\right) \nonumber \\ 
& {} + \frac{p_6p_7(p_2+p_3)^2}{(p_6+p_7)}(N-1)(N-2)(1-p_6-p_7)^{N-3} F\left([1, 1, 3-N], [2, 2], -\frac{p_6+p_7}{1-p_6-p_7}\right) \nonumber \\ 
& {} + \frac{p_6p_7}{(p_6 + p_7)^2}(p_2+p_3 + p_Z) + \frac{p_7^2}{(p_6+p_7)^2}p_{Z}(1-p_{Z})
\end{align}

\sbsbsection{The Covariances}
For the covariances where separate generating variables are used
\begin{align}
E[M_{11}M_{10}] = \frac{1}{N}\int\mathrm{d}w \frac{ax}{w}\frac{\partial^2}{\partial a\partial w}\left[\frac{bzp_7(p_2u+p_3uv+p_6y+p_7yz)}{(p_6 + p_7z)}\right] T_{N-1}\Bigg|_{\substack{a=b=1 \cr s=t=1 \cr u=v=1 \cr w=x=1 \cr y=z=1}}
\end{align}
it is easy to see that the operators act on $T_{N-1}$ rather than the prefactor, so we repeat the calculation for the mean with $N$ replaced by $(N-1)$ to obtain
\begin{align}
C[M_{11},M_{10}] = -\frac{1}{N}E[M_{11}]E[M_{10}] \, , & & C[M_{01},M_{10}] = -\frac{1}{N}E[M_{01}]E[M_{10}] \nonumber \\ 
C[M_{11},M_{00}] = -\frac{1}{N}E[M_{11}]E[M_{00}] \, , & & C[M_{01},M_{00}] = -\frac{1}{N}E[M_{01}]E[M_{00}]
\end{align}
The more complicated cases are where the generating variables reoccur
\begin{align}
E[M_{11}M_{01}] = \frac{1}{N}\int\mathrm{d}u \frac{bv}{u}\frac{\partial^2}{\partial b\partial v}\left[\frac{bzp_7u(p_2+p_3v)}{(p_6 + p_7z)} + bzp_7y\right] T_{N-1}\Bigg|_{\substack{a=b=1 \cr s=t=1 \cr u=v=1 \cr w=x=1 \cr y=z=1}}
\end{align}
For the term linear in $u$ we first integrate then differentiate wrt $v$ while for the other term we first differentiate then integrate to give
\begin{align}
E[M_{11}M_{01}] &= \frac{1}{N}\frac{\partial}{\partial b} \left[\frac{bp_3p_7}{(p_6 + p_7)} + \frac{bp_3p_7}{(p_2 + p_3)}\right] T_{N-1}\Bigg|_{\substack{a=b=1 \cr s=t=1 \cr u=v=1 \cr w=x=1 \cr y=z=1}} \nonumber \\
& = E[M_{11}]E[M_{01}] + \frac{p_3p_7}{N(p_2+p_3)(p_6+p_7)}p_{Z}(1-p_Z)
\end{align}
and
\begin{align}
C[M_{11},M_{01}] &= \frac{1}{N}\frac{p_3p_7}{(p_2+p_3)(p_6+p_7)}p_{Z}(1-p_Z) \nonumber \\
C[M_{01},M_{00}] &= \frac{1}{N}\frac{p_1p_5}{(p_0+p_1)(p_4+p_5)}p_{Z}(1-p_Z)
\end{align}

\sbsbsection{The Variance}
Since the terms from the covariances simplify, the complete variance is
\begin{align}
& V \cdot N = \frac{p_6p_7(p_2+p_3)}{(p_6+p_7)}(N-1)(1-p_6-p_7)^{N-2} F\left([1, 1, 2-N], [2, 2], -\frac{p_6+p_7}{1-p_6-p_7}\right) \nonumber \\ 
& + \frac{p_6p_7(p_2+p_3)^2}{(p_6+p_7)}(N-1)(N-2)(1-p_6-p_7)^{N-3} F\left([1, 1, 3-N], [2, 2], -\frac{p_6+p_7}{1-p_6-p_7}\right) \nonumber \\ 
& + \frac{p_4p_5(p_0+p_1)}{(p_4+p_5)}(N-1)(1-p_4-p_5)^{N-2} F\left([1, 1, 2-N], [2, 2], -\frac{p_4+p_5}{1-p_4-p_5}\right) \nonumber \\ 
& + \frac{p_4p_5(p_0+p_1)^2}{(p_4+p_5)}(N-1)(N-2)(1-p_4-p_5)^{N-3} F\left([1, 1, 3-N], [2, 2], -\frac{p_4+p_5}{1-p_4-p_5}\right) \nonumber \\ 
& + \frac{p_2p_3(p_6+p_7)}{(p_2+p_3)}(N-1)(1-p_2-p_3)^{N-2} F\left([1, 1, 2-N], [2, 2], -\frac{p_2+p_3}{1-p_2-p_3}\right) \nonumber \\ 
& + \frac{p_2p_3(p_6+p_7)^2}{(p_2+p_3)}(N-1)(N-2)(1-p_2-p_3)^{N-3} F\left([1, 1, 3-N], [2, 2], -\frac{p_2+p_3}{1-p_2-p_3}\right) \nonumber \\ 
& + \frac{p_0p_1(p_4+p_5)}{(p_0+p_1)}(N-1)(1-p_0-p_1)^{N-2} F\left([1, 1, 2-N], [2, 2], -\frac{p_0+p_1}{1-p_0-p_1}\right) \nonumber \\ 
& + \frac{p_0p_1(p_4+p_5)^2}{(p_0+p_1)}(N-1)(N-2)(1-p_0-p_1)^{N-3} F\left([1, 1, 3-N], [2, 2], -\frac{p_0+p_1}{1-p_0-p_1}\right) \nonumber \\ 
& + \frac{p_6p_7}{(p_6 + p_7)^2}(p_2+p_3 + p_Z) + \frac{p_4p_5}{(p_4 + p_5)^2}(p_0+p_1 + 1 - p_Z)  \nonumber \\
& + \frac{p_2p_3}{(p_2 + p_3)^2}(p_6+p_7 + p_Z) + \frac{p_0p_1}{(p_0 + p_1)^2}(p_4+p_5 + 1 - p_Z)  \nonumber \\
& + \left[\frac{p_7}{(p_6+p_7)} - \frac{p_5}{(p_4+p_5)} - \frac{p_3}{(p_2+p_3)} + \frac{p_1}{(p_0+p_1)}\right]^{2}p_{Z}(1-p_{Z}) 
\end{align}
We note that the hypergeometric functions can be written solely in terms of $p_X$ and $p_Z$ so that the variance is actually quadratic in $p_{Y,i}$.

\section{Asymptotics of the hypergeometric functions} \label{asymphyper}

We utilise the following asymptotic expansions of our hypergeometric functions:
\begin{align}
N^{2}z^{2}(1-z)^{N-1}F\left([1, 1, 1-N], [2, 2], -\frac{z}{1-z}\right) = 1 + \frac{(1-z)}{Nz} + \ldots
\end{align}
and
\begin{align}
(N-1)z^{2}(1-z)^{N-2}F\left([1, 1, 2-N], [2, 2], -\frac{z}{1-z}\right) & = \frac{1}{N} + \ldots \nonumber \\
(N-1)(N-2)z^{2}(1-z)^{N-3}F\left([1, 1, 3-N], [2, 2], -\frac{z}{1-z}\right) & = 1 + \frac{1}{Nz} + \ldots
\end{align}

\end{document}